\documentclass[11pt]{article}

\usepackage{amsmath}    
\usepackage{graphicx}   
\usepackage{verbatim}   
\usepackage{color}      
\usepackage{subfigure}  
\usepackage{amssymb,tikz}
\usepackage{soul}

\usepackage[bookmarks, bookmarksopen, bookmarksdepth=3, colorlinks=true,
linkcolor=black, anchorcolor=black, citecolor=black, filecolor=black,
menucolor=black, runcolor=black, urlcolor=black, pdfencoding=unicode
]{hyperref}

\newtheorem{theorem}{Theorem}[section]
\newtheorem{proposition}[theorem]{Proposition}

\newtheorem{lemma}[theorem]{Lemma}

\newtheorem{problem}[theorem]{Problem}

\newcommand{\proof}{\noindent{\bf Proof.\ }}
\newcommand{\qed}{\hfill $\square$\medskip}

\newcommand{\R}{\mathbb R}
\newcommand{\N}{\mathbb N}

\newcommand{\V}{\mathcal V}

\renewcommand{\U}{\mathcal U}
\newcommand{\SSS}{\mathcal S}
\renewcommand{\gg}{\gamma_{g}}
\newcommand{\gtg}{{\gamma_{g}^{t}}}

\let\d\relax
\DeclareMathOperator {\d} {d}

\begin{document}

\tikzstyle{every node}=[circle, draw, fill=black!10,
                        inner sep=0pt, minimum width=4pt]

\title{Effect of predomination and vertex removal on the game total domination number of a graph}
 
 \author{
 	Vesna Ir\v si\v c $^{a, b}$
 }
 
 \date{\today}
 
 \maketitle
 \begin{center}
 	$^a$ Institute of Mathematics, Physics and Mechanics, Ljubljana, Slovenia\\
 	{\tt vesna.irsic@fmf.uni-lj.si }
 	\medskip
 
 	$^b$ Faculty of Mathematics and Physics, University of Ljubljana, Slovenia\\
 		
 \end{center}
 
 \begin{abstract}
 The game total domination number, $\gtg$, was introduced by Henning et al.\ in 2015. In this paper we study the effect of vertex predomination on the game total domination number. We prove that $\gtg(G|v) \geq \gtg(G) - 2$ holds for all vertices $v$ of a graph $G$ and present infinite families attaining the equality. To achieve this, some new variations of the total domination game are introduced. The effect of vertex removal is also studied. We show that $\gtg(G) \leq \gtg(G-v) + 4$ and $\gtg'(G) \leq \gtg'(G-v) + 4$. 
 \end{abstract}
 
 
 \medskip
 \noindent{\bf AMS Subj.\ Class.:} 05C57, 05C69

\section{Introduction}
\label{sec:intro}
 
The domination game was introduced in 2010 by Bre\v{s}ar et al.~\cite{dom} as a game played by two players, Dominator and Staller, on the graph $G$. They alternate taking turns for as long as possible and on each turn one chooses such a vertex in $G$ that dominates at least one not yet dominated vertex. Recall that a vertex dominates itself and its neighbors. Dominator tries to minimize and Staller tries to maximize the number of moves. The total number of selected vertices is called \emph{the game domination number}, $\gg(G)$, if Dominator starts the game (\emph{D-game}) or \emph{the Staller-start game domination number}, $\gg'(G)$, if Staller makes the first move (\emph{S-game}). If Staller is allowed to pass one move, the game is called the \emph{Staller-pass game} and the number of moves made is $\gg^{sp} (G)$ or $\gg'^{sp} (G)$. Analogously we define the \emph{Dominator-pass game}. Graphs on which an optimal domination game yields a minimum dominating set of the graph have been introduced and studied in~\cite{D-trivial}. See~\cite{XLK-2018}, for results about graphs with maximal possible game domination number, i.e.\ $2 \gamma (G) - 1$.

One can also consider a domination game played on the graph $G$ where some vertices are already considered dominated~\cite{dom}. If $S \subseteq V(G)$ are already dominated, then the resulting game domination number is denoted by $\gg(G|S)$ or $\gg'(G|S)$. If $S = \{v\}$, then we write $\gg(G|v)$ or $\gg'(G|v)$. An important property arising from this definition is the Continuation Principle~\cite{extremal} which states the following. For a graph $G$ and the sets $A, B \subseteq V(G)$, $B \subseteq A$, it holds that $\gg(G|A) \leq \gg(G|B)$ and $\gg'(G|A) \leq \gg'(G|B)$.

The total version of the domination game was introduced in~\cite{totDom} and is defined analogously, except that on each turn only a vertex which totally dominates at least one not yet totally dominated vertex can be played. Recall that a vertex totally dominates only its neighbors and not itself. \emph{The game total domination number}, $\gtg(G)$, is the number of vertices chosen if Dominator starts the game on $G$. If Staller plays first, then \emph{the Staller-start game total domination number} is denoted by $\gtg'(G)$. As in the domination game, \emph{Staller- and Dominator-pass games}, as well as games on partially dominated graphs, can also be considered. See~\cite{trees}, for results about the game total domination number on trees.

An important property of the game total domination number is the Total Continuation Principle~\cite{totDom}, which states that for a graph $G$ and the sets $A, B \subseteq V(G)$, $B \subseteq A$, it holds that $\gtg(G|A) \leq \gtg(G|B)$ and $\gtg'(G|A) \leq \gtg'(G|B)$.  Another fundamental property, which is also parallel to the ordinary domination game, is the following~\cite{totDom}. For any graph $G$, we have $|\gtg(G) - \gtg'(G)| \leq 1$.

The game total domination number is non-trivial even on paths and trees~\cite{totDomCP} and is log-complete in PSPACE~\cite{pspace}. The $\frac{3}{4}$-conjecture, stating that for all graphs $G$ it holds $\gtg(G) \leq \frac{3}{4} |V(G)|$, was posed in~\cite{3/4_1} and further studied in~\cite{3/4_2, 3/4_3}.

\emph{Game domination critical} (or \emph{$\gg$-critical}) graphs have been introduced in~\cite{DomCritical} as graphs $G$ for which $\gg(G) > \gg(G|v)$ for all $v \in V(G)$. Among other results from~\cite{DomCritical}, we recall that for any vertex $u \in V(G)$ it holds $\gg(G|u) \geq \gg(G) - 2$ and that there exist a graph attaining the equality. Analogously, \emph{total domination game critical} (or \emph{$\gtg$-critical}) graphs were introduced in~\cite{totDomCritical} as graphs $G$ for which $\gtg(G) > \gtg(G|v)$ for all $v \in V(G)$. Infinite families of $\gtg$-critical circular and M\"{o}bius ladders were presented in~\cite{ladders}.

In Section~\ref{sec:sth} we state that $\gtg(G|v) \geq \gtg(G) - 2$ for every $v \in V(G)$ and observe that neither of the total domination game critical graphs studied in~\cite{totDomCritical, ladders} attains the equality. Hence, in Sections~\ref{sec:sth} and~\ref{sec:generalisation} we present infinite families of graphs attaining it. We apply the introduced techniques to another family of graphs as well.

A similar concept, the effect of edge or vertex removal on the game domination number has been studied in~\cite{EdgeVertexRemoval, G-e, book}. It holds that $\gtg(G-v)$ cannot be bounded from above by $\gtg(G)$, but that $\gtg(G) - \gtg(G-v) \leq 2$. Examples of graphs attaining $\gtg(G) - \gtg(G-v) \in \{0,1,2\}$ have also been presented. Analogous results hold for the Staller-start game. In Section~\ref{sec:vertex} we present similar results for the game total domination number.

\section{Graphs with the property $\gtg(G|v) = \gtg(G) - 2$}
\label{sec:sth}

As mentioned, the game total domination critical graphs were introduced in~\cite{totDomCritical}, where also critical cycles and paths were characterized, as well as $2$- and $3$-$\gtg$-critical graphs. Domination game critical graphs were introduced and studied in~\cite{DomCritical}, where we find the following property: $\gamma_g(G | u) \geq \gamma_g (G) - 2$ for every vertex $u$. With the same reasoning, we derive an analogous result for the game total domination number. To be self contained, we rephrase the proof here without detail.

\begin{lemma}
	\label{lem:-2}
	For every $v \in V(G)$, it holds $\gtg(G|v) \geq \gtg(G) - 2$.
\end{lemma}

\proof
The players play the real game on $G$ while Dominator imagines a Staller-pass game on $G|u$. He ensures that all the vertices that are totally dominated in the real game are also totally dominated in the imagined game. Dominator plays optimally in the imagined game and copies each of his moves into the real game. Every move of Staller in the real game is copied into the imagined game. If this move is not legal in the imagined game, the only new totally dominated vertex in the real game in $u$. In this case, Staller simply skips her move in the imagined game. 

Let $p$ denote the number of moves played in the real game and $q$ denote the number of moves played in the imagined game. As Staller plays optimally in the real game, it holds $\gtg(G) \leq p$. As in the imagined game one move might be skipped, we have $p \leq q + 1$. As Dominator plays optimally on $G | u$, it holds $q \leq \gtg^{sp} (G|u)$. Similarly as in~\cite{union}, we can also derive $\gtg^{sp}(G|u) \leq \gtg(G|u) + 1$. Combining these inequalities yields $\gtg(G) \leq \gtg(G|u) + 2$.
\qed

For all game total domination critical graphs in~\cite{totDomCritical, ladders} it holds that for each vertex $v$ we have $\gtg(G|v) = \gtg(G) - 1$. A natural question is whether there exist graphs $G$ with $\gtg(G|v) = \gtg(G) - 2$ for some vertex $v \in V(G)$, i.e.\ graphs attaining the equality in Lemma~\ref{lem:-2} for some vertex $u$. It turns out that the answer is positive. To show this, we consider the following family of graphs.

\begin{figure}[!ht]
 \begin{center}
 \begin{tikzpicture}[thick,scale=0.5]
 
 \pgfmathtruncatemacro{\N}{14}
 \pgfmathtruncatemacro{\M}{4}
 \pgfmathtruncatemacro{\r}{2}
 \pgfmathtruncatemacro{\R}{4}
 \begin{scope}
 \foreach \x in {1,...,\N}
 \node[label={-80 + \x*360/\N:$u_{\x}$}] (\x) at (-80 + \x*360/\N:\R cm) {};
 \foreach \x [remember=\x as \lastx (initially 1)] in {1,...,\N,1}
 \path (\x) edge (\lastx);
 \end{scope}
 
 \node[label=below: {$w_1$}] (w1) at (10, -2) {};
 \node[label=above: {$w_2$}] (w2) at (10, 2) {};
 \node[label=below: {$v_1$}] (v1) at (14, -2) {};
 \node[label=above: {$v_2$}] (v2) at (14, 2) {};
 
 \path (1) edge (w1);
 \path (5) edge (w1);
 \path (1) edge (w2);
 \path (5) edge (w2);
 
 \path (v1) edge (w1);
 \path (v2) edge (w1);
 \path (v1) edge (w2);
 \path (v2) edge (w2);
 \path (w2) edge (w1);
 \path (v2) edge (v1);

 \end{tikzpicture}
 \caption{A graph $G_{14, 4}$.}
 \label{fig:primer}
 \end{center}
 \end{figure}
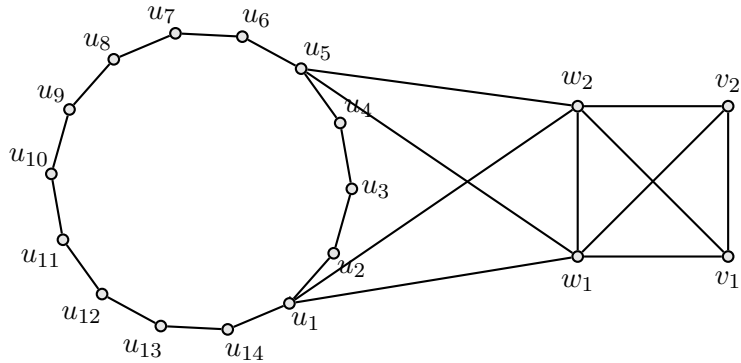

Let $G_{n, m}$, $n \geq 8, m \geq 4$, $n \equiv 2 \mod 6$, be a graph consisting of a cycle $C_n$ with vertices $\U = \{u_1, \ldots, u_n\}$ (and naturally defined edges), a complete graph $K_m$ with vertices $w_1, w_2$ and $\V = \{ v_1, \ldots, v_{m-2}\}$ (and naturally defined edges), and additional edges $\{u_1, w_1\}, \{u_1, w_2\}, \{u_5, w_1\}, \{u_5, w_2\}$ (cf.\ Fig.~\ref{fig:primer}). Denote by $H_m$ a subgraph induced on vertices $\{u_1, u_5, w_1, w_2\} \cup \V$. Observe that $G_{n,m}$ is $2$-connected. We shall prove the following.

\begin{theorem}
	\label{thm:minus2}
	For $n \geq 8, m \geq 4$, $n \equiv 2 \mod 6$, it holds that $\gtg(G_{n, m}) = \frac{2n-1}{3} + 2$ and $\gtg(G_{n, m}|w_1) = \frac{2n-1}{3}$.
\end{theorem}

During the proof we will consider several different variations of the game played on the cycle. We will study them seperately.

\subsection{Preliminaries}
\label{sec:preliminaries}

Recall from~\cite{totDomCP, totDomCritical} that for $n \equiv 2 \mod 6$ it holds 
$$\gtg(C_n) = \gtg(C_n|v) = \frac{2n-1}{3}
\quad \text{and} \quad
\gtg'(C_n) = \gtg'(C_n|v) = \frac{2n-1}{3} - 1$$
for every vertex $v \in V(C_n)$. In some of the following proofs we will use the strategy for Staller from~\cite{totDomCP}, thus we rephrase it here. First recall that a \emph{run} on a partially dominated cycle is a maximal sequence of (at least two) consecutive totally dominated vertices. An \emph{anti-run} is a maximal sequence of (at least two) consecutive not totally dominated vertices. Let $A$ denote the set of already totally dominated vertices on a cycle $C_n$. Suppose $A$ is neither empty nor $V(C_n)$. If $A$ contains a run or an anti-run, Staller can play on its extremity and totally dominate only one new vertex. If $A$ contains no runs and no anti-runs, then $(A, A^C)$ must be a bipartition of the cycle. In this case, Staller cannot totally dominate only one new vertex in the next move. We call this strategy $S_1$. 

Recall also the strategy for Dominator from~\cite{totDomCP}. Let $U$ denote the set of unplayable vertices. Clearly, Staller adds at least one vertex to $U$ on each move. Dominator's strategy is to ensure that after her move and his reply, at least three vertices are added to $U$. Say Staller plays on $v_1$. She totally dominates at least one new vertex, say $v_2$. Now label the other vertices cyclically as $v_3, \ldots, v_n$. Since $v_2$ was not yet totally dominated before Staller played $v_1$, vertex $v_3$ was not played yet. If $v_5$ was also not played yet, then $v_4$ is not yet totally dominated and Dominator can play on $v_5$. Thus, he adds $v_1, v_3, v_5$ to $U$. But if $v_5$ was played before, then Staller's move makes both $v_1$ and $v_3$ unplayable. Thus, Dominator can reply anywhere to add another vertex to $U$. We call this strategy $D_1$. 

Notice also that
$$\gtg(H_m) = \gtg'(H_m) = \gtg'(H_m|v) = 2$$
for every $v \in V(H_m)$ and
$$\gtg(H_m|v) = \begin{cases}
1; & v \in \{w_1, w_2\},\\
2; & \text{otherwise}.
\end{cases} $$

Now we prove a series of lemmas.

\begin{lemma}
	\label{lem:Cycle15}
	If $n \equiv 2 \mod 6$, then $\gtg(C_n|\{u_1, u_5\}) = \frac{2n-1}{3}$.
\end{lemma}

\proof
It follows from the total continuation principle that $\gtg(C_n|\{u_1, u_5\}) \leq \gtg(C_n) = \frac{2n-1}{3}$. Thus we only need to prove $\gtg(C_n|\{u_1, u_5\}) \geq \frac{2n-1}{3}$, which can be done by finding a suitable strategy for Staller.

Let $n = 6 q + 2$ for some positive integer $q$. Denote $x_k = u_{2k-1}$, $y_k = u_{2k}$ for $k \in [\frac{n}{2}]$ and $X = \{x_1, \ldots, x_k\}$, $Y = \{y_1, \ldots, y_k\}$. Clearly, $(X, Y)$ is a bipartition of the cycle $C_n$. So by playing in $X$, a player can totally dominate only vertices in $Y$, and vice versa. As $u_1, u_5 \in X$ are already predominated, the players need to dominate only $3q-1$ vertices in $X$, and $3q+1$ vertices in $Y$. 

Observe that if Dominator plays on $X$ and after his move $Y$ is not yet totally dominated, then there exists $y_i, y_{i+1}$ such that one is already totally dominated and the other is not. Then Staller can play their common neighbor $x_{i+1} \in X$ and hence she totally dominates only one new vertex in such a move. A similar observation holds if we switch $X$ and $Y$. We now distinguish two cases.

\begin{enumerate}
	\item Dominator is the first player to play on $X$.\\
	 Staller's strategy is to reply on $X$ whenever Dominator plays on $X$, and on $Y$, whenever he plays on $Y$. Suppose Dominator makes in total $\ell$ moves on $X$. So Staller makes $\ell$ or $\ell-1$ moves on $X$ and it follows from the above observation, that she can totally dominate only one new vertex on each move. 
	
	Suppose Dominator totally dominates two new vertices on each move on $X$. Then together both players dominate $2\ell + \ell = 3\ell$ or $2\ell + (\ell-1) = 3\ell-1$ vertices in $Y$. But as $|Y| \equiv 1 \mod 3$, this is not possible. Hence, Dominator has to dominate only one new vertex on at least one move on $X$. 
	
	If Staller plays first on $Y$, she plays $y_1$ (i.e.\ $u_2)$, so she only totally dominates one new vertex. In all other moves, or if Dominator starts playing on $Y$, it is clear that Staller can totally dominate only one new vertex on each move. 
	
	Let $m$ denote the number of moves played. If $m$ is odd, we have 
	$$n \leq 2 + \frac{m-1}{2} \cdot 1 + \frac{m+1}{2} \cdot 2 - 1 = \frac{3m+3}{2} \,,$$
	so $$m \geq \left\lceil \frac{2n-3}{3} \right\rceil = \frac{2n-1}{3}\,.$$
	
	If $m$ is even, we have 
	$$n \leq 2 + \frac{m}{2} \cdot 1 + \frac{m}{2} \cdot 2 - 1 = \frac{3m+2}{2}\,,$$
	thus $$m \geq \left\lceil \frac{2n-2}{3} \right\rceil = \frac{2n-1}{3}\,.$$
	
	\item Staller is the first player to play in $X$.\\
	 Thus Dominator started on $Y$ and then Staller replies on $X$ (resp.\ $Y$) if Dominator plays in $X$ (resp.\ $Y$). But she follows additional strategy to avoid playing on $\{y_1, y_2\}$ (i.e.\ $\{u_2, u_4\}$) for as long as possible. Notice that she can always play elsewhere unless all other vertices in $X-\{x_2\}$ are already totally dominated, as $x_2$ is the common neighbor of $y_1, y_2$, and $x_1, x_3$ are already totally dominated. But as Staller should play first on $X$, Dominator has to play last on $Y$ and as $x_2$ will not be dominated by Staller (due to her strategy), it is clear that Dominator will totally dominate only one new vertex on at least one of his moves on $Y$.   
	
	As above it is clear that Staller can totally dominate only one new vertex on each of her moves, except when she plays first on $X$. 
	
	Denote by $\ell$ the number of moves Dominator makes on $Y$. Suppose he totally dominates just one new vertex on only one of his moves on $Y$. Because he should also play last on $Y$, Staller makes $\ell-1$ moves on $Y$. So together they totally dominate $2\ell - 1 + \ell-1 = 3\ell-2$ vertices in $X$, which is a contradiction with the fact that $3q-1$ vertices in $X$ need to be totally dominated.
	
	Therefore, Dominator totally dominates just one new vertex on at least two of his moves. Let again $m$ denote the number of moves played. If $m$ is odd, we have 
	$$n \leq 2 + \frac{m-1}{2} \cdot 1 + 1 + \frac{m+1}{2} \cdot 2 - 2 = \frac{3m+3}{2}\,,$$
	so $$m \geq \left\lceil \frac{2n-3}{3} \right\rceil = \frac{2n-1}{3}\,.$$
	
	If $m$ is even, we have 
	$$n \leq 2 + \frac{m}{2} \cdot 1 + 1 + \frac{m}{2} \cdot 2 - 2 = \frac{3m+2}{2}\,,$$
	thus $$m \geq \left\lceil \frac{2n-2}{3} \right\rceil = \frac{2n-1}{3}\,.$$
\end{enumerate}
Hence also $\gtg(C_n|\{u_1, u_5\}) \geq \frac{2n-1}{3}$.
\qed

Consider the following variation of the total domination game---Staller plays twice and only then the players start to alternate moves. So we have moves $s_1, s_2, d_1, s_3, \ldots$ The number of moves in such a game is denoted by $\gtg''(G)$. 

\begin{lemma}
	\label{lem:ss}
	If $n\geq 3$ is a positive integer, then $\gtg''(C_n) = \gtg(C_n)$. 
\end{lemma}

\proof
The first move $s_1$ of Staller totally dominates two new vertices. After that the players alternate moves and both play optimally. But as the cycle is vertex-transitive, the first move $s_1$ of Staller can be considered as an optimal first move of Dominator in the usual D-game. Hence, $\gtg''(C_n) = \gtg(C_n)$.
\qed

Consider another variation of the total domination D-game---both players alternate moves normally, but after $m$, $0 \leq m \leq \gtg(G)$, moves of the game  vertices $x_1, \ldots, x_k$ become totally dominated (for free). The number of moves in such a game is denoted by $\gtg(G|^{m}\{x_1, \ldots, x_k\})$. Notice that $\gtg(G|^{0}\{x_1, \ldots, x_k\}) = \gtg(G|\{x_1, \ldots, x_k\})$ and $\gtg(G|^{\gtg(G)}\{x_1, \ldots, x_k\}) = \gtg(G)$. We point out that neither of the players is aware in advance of the parameter $m$ and the set $\{x_1, \ldots, x_k\}$.

\begin{lemma}
	\label{lem:cycle15}
	If $n \equiv 2 \mod 6$, $n \geq 8$, then $\gtg(C_n|^{m}\{u_1, u_5\}) \geq \frac{2n-1}{3}$ for any $m$, $0 \leq m \leq \gtg(C_n)$.
\end{lemma}

\proof
The real game is played simultaneously with an imagined total domination D-game on $C_n|\{u_1, u_5\}$ which is imagined by Staller. We will describe Staller's strategy that will ensure that the set of totally dominated vertices in the real game is a subset of the totally dominated vertices in the imagined game. 

Each move of Dominator is copied to imagined game. If it is not playable, any other legal move is played there. Staller replies optimally in the imagined game and copies her move to the real game. This is always legal, even after $u_1$ and $u_5$ become predominated in the real game. Hence by Lemma~\ref{lem:Cycle15}, it holds $\gtg(C_n|^{m}\{u_1, u_5\}) \geq \gtg(C_n|\{u_1, u_5\}) =  \frac{2n-1}{3}$ for any $m$.
\qed

Another variation we shall introduce is the following. Consider a total domination D-game on a graph $G$, where Staller has to pass one move, and after that also Dominator has to pass one move. The number of moves in such a game is denoted by $\gtg^{sdp (k,\ell)}(G)$ if Staller passes exactly after $k$-th move of the game and Dominator passes exactly after the $\ell$-th move, $k \leq \ell$. These two passes are not counted as moves.

\begin{lemma}
	\label{lem:sdp}
	If $n \equiv 2 \mod 6$, $n \geq 8$, then $\gtg^{sdp (k,\ell)}(C_n) \geq \frac{2n-1}{3}$ for any $0 \leq k \leq \ell \leq \gtg(C_n)$.
\end{lemma}

\proof
Following the same strategy as presented in~\cite{totDomCP}, Staller can totally dominate only one vertex on each of her moves, except maybe once. Let $m$ be the number of moves played.

If $m$ in even, we have $n \leq \frac{m}{2} \cdot 1 + 1 + \frac{m}{2} \cdot 2 = \frac{3m+2}{2}$, and thus $m \geq \frac{2n-1}{3}$.

If $m$ is odd, we have $n \leq \frac{m-1}{2} \cdot 1 + 1 + \frac{m+1}{2}\cdot 2 = \frac{3m+3}{2}$, and thus $m \geq \frac{2n-1}{3}$.
\qed

The combination of the above variations is such that Staller has to pass one move, and then Dominator has to pass one move, but when Dominator passes, some vertices $x_1, \ldots, x_k$ become predominated. The number of moves in such a game is denoted by $\gtg^{sdp (m,\ell)}(G|^{\ell}\{x_1, \ldots, x_k\})$ if Staller passes after the $m$-th move of the game and Dominator passes after the $\ell$-th move, $m \leq \ell$. Combining Lemma~\ref{lem:sdp} with the imagination strategy used in the proof of Lemma~\ref{lem:cycle15}, yields the following lemma.

\begin{lemma}
	\label{lem:sdp15}
	If $n \equiv 2 \mod 6$, $n \geq 8$, then $\gtg^{sdp(m,\ell)}(C_n|^{\ell}\{u_1, u_5\}) \geq \frac{2n-1}{3}$.
\end{lemma}

Consider now our final variation. Consider a D-game on $G$ where Staller has to pass exactly twice. Each pass appears after Dominator plays the vertex $u'$ or $u''$ for some $u', u'' \in V(G)$. Additionally, we allow Dominator to play his first move on $\{u', u''\}$ even if this dominates no new vertices on $G$. But his second move on $\{u', u''\}$ has to totally dominate some new vertex in $G$. The number of moves in such a game is denoted by $\gtg^{ssp, u', u''}(G)$. In the case when $G$ is a cycle $C_n$, we can write $\gtg^{ssp, u', u''}(C_n) = \gtg^{ssp, \d(u', u'')}(C_n)$, where $\d(u', u'')$ denotes the distance between the special vertices, $u'$ and $u''$. In our case, $\d(u_1, u_5) = 4$.

\begin{lemma}
	\label{lem:ssp}
	If $n \equiv 2 \mod 6$, $n \geq 8$, then $\gtg^{ssp,4}(C_n) \geq \frac{2n-1}{3}$.
\end{lemma}

\proof
Staller again follows her strategy from~\cite{totDomCP}, which ensures that except perhaps on one move, she can totally dominate just one new vertex. Let $m$ denote the total number of moves played. Distinguish two cases.

\begin{enumerate}
	\item Staller can always totally dominate just one new vertex.\\
	If $m$ in even, we have $n \leq \frac{m-2}{2} \cdot 1 + \frac{m+2}{2} \cdot 2 = \frac{3m+2}{2}$, and thus $m \geq \frac{2n-1}{3}$.
	
	If $m$ is odd, we have $n \leq \frac{m-3}{2} \cdot 1  + \frac{m+3}{2} \cdot 2 = \frac{3m+3}{2}$, and thus $m \geq \frac{2n-1}{3}$.
	
	\item Staller dominates two new vertices on some move.\\
	We know this happens exactly when it is Staller's turn when $(A, A^C)$ forms a bipartition of the cycle, where $A$ denotes already totally dominated vertices. If it is not a bipartition, Staller can totally dominate only one new vertex (due to the strategy $S_1$).
	
	If Dominator totally dominates just one (or zero) new vertex on some move, calculations similar to the above show $m \geq \frac{2n-1}{3}$. From now on suppose Dominator totally dominates exactly two new vertices on each move.
	
	If $(A, A^C)$ is a bipartition before $u_1$ and  $u_5$ were played, it follows $u_1 \in A$. Hence Staller can play $u_3$. But than Dominator cannot play $u_1$ or $u_5$ and still dominate two vertices on his move. 
	
	If $(A, A^C)$ is a bipartition after $u_1$ and $u_5$ were played, it follows that $u_1 \in A^C$. So both Staller's passes appeared before this situation and also it has to be Staller's turn now. Denote by $d$ the number of moves Dominator made before $(A, A^C)$ is a bipartition. Then Staller made only $d-3$ moves up to this point. So together they totally dominated $2d + (d-3) = 3 (d-1)$ vertices, which is a contradiction, as $|A| \equiv 1 \mod 3$. 
\end{enumerate}
Thus in all cases we have $m \geq \frac{2n-1}{3}$.
\qed

%

\subsection{Proof of Theorem~\ref{thm:minus2}}
\label{sec:proof}

First we describe a strategy for Dominator to show that $\gtg(G_{n, m}) \leq \frac{2n-1}{3} + 2$. Dominator starts by playing $v_1$. If Staller replies on $K_m$, the whole $K_m$ is totally dominated after her move. So what remains is a normal game on $C_n$ or $C_n|\{u_1, u_5\}$, thus Dominator just plays according to his optimal strategy there. In this case the total number of moves is $\frac{2n-1}{3} + 2$. 

But if Staller replies on $C_n$, Dominator's second move is $v_2$. After his move, the whole $K_m$ is totally dominated and we get a $\gtg''$-game on $C_n$. It follows from Lemma~\ref{lem:ss} that in this case $\frac{2n-1}{3} + 2$ moves are made. Hence, $\gtg(G_{n, m}) \leq \frac{2n-1}{3} + 2$.

Next we describe a strategy for Staller to show that $\gtg(G_{n, m}) \geq \frac{2n-1}{3} + 2$. If Dominator plays on $\{u_1, u_5\}$ or on  $\{w_1, w_2\}$, Staller replies on $\V$. If Dominator plays on $\V$, Staller replies on $\{w_1, w_2\}$. If the prescribed reply is not legal, Staller plays any legal move. Note that this can only happen if the whole $K_m$ is already totally dominated. 

If Dominator plays on $C_n - \{u_1, u_5\}$, Staller replies on such optimal move on $C_n$, that she does not play on $\{u_1, u_5\}$. Why is this possible whenever $C_n$ is not yet totally dominated? Let $\SSS$ denote the set of all Staller's optimal moves on $C_n$. Due to strategy $S_1$, these are the vertices which have only one not yet totally dominated neighbor, or half of all vertices in the case when $(A, A^C)$ forms a bipartition. If $\SSS \cap \{u_1, u_5\} \neq \emptyset$, Staller can make an optimal move on $C_n - \{u_1, u_5\}$. Else, $\SSS \subseteq \{u_1, u_5\}$ and only one new vertex will be totally dominated. Without loss of generality, $u_1 \in \SSS$. Suppose $u_2$ is already totally dominated and $u_n$ is not. But as $u_{n-1} \notin \SSS$, $u_{n-2}$ is also not totally dominated, otherwise Staller could play $u_{n-1}$. Repeating this reasoning results in the fact that $u_2$ is not totally dominated, which is a contradiction. Supposing that $u_2$ is not totally dominated and $u_n$ is, yields a similar contradiction. Hence, Staller can play on $C_n - \{u_1, u_5\}$, as long as $C_n$ is not yet totally dominated. If it is, she plays on $\V$.

This strategy ensures that at least two moves are played on $K_m$. Thus we only need to prove that at least $\frac{2n-1}{3}$ moves are played on $C_n$. Consider now only moves played on the cycle. If Dominator does not start on it, after Staller's reply the whole $K_m$ is totally dominated and all remaining legal moves are on $C_n$. Hence, this is always a D-game. We distinguish two cases.

\begin{enumerate}
	\item Dominator is the first player to play on $K_m$.\\
	After Staller's reply, the whole $K_m$ is totally dominated so exactly two moves are played on $K_m$. Thus we have a normal game on $C_n$, where possibly vertices $\{u_1, u_5\}$ become predominated at some point of the game. The number of moves on $C_n$ is therefore at least $\gtg(C_n)$ or $\gtg(C_n|^{\ell}\{u_1, u_5\})$ for some $\ell$. It follows from Lemma~\ref{lem:cycle15}, that in either case at least $\frac{2n-1}{3}$ moves are played on $C_n$.
	
	\item Staller is the first player to play on $K_m$. 
	\begin{enumerate}
		\item If the whole $C_n$ is totally dominated before Staller first plays on $K_m$ (specifically, on $\V$), a normal D-game was played on $C_n$, thus at least $\frac{2n-1}{3}$ moves were made on the cycle.
		
		\item Otherwise Dominator played on $\{u_1, u_5\}$ right before Staller's move on $K_m$. Dominator can reply on:
		\begin{enumerate}
			\item $\V$: with this move he completes the total domination of $K_m$ without any vertex of $C_n$ becoming predominated. This results on a normal $\gtg$-game on the cycle.
			
			\item $\{w_1, w_2\}$: similar as above, but vertices $\{u_1, u_5\}$ become predominated, so at least $\gtg(C_n|^{\ell}\{u_1, u_5\})$ moves are played on the cycle (for some $\ell$).
			
			\item $C_n$: if Dominator plays on $\{u_1, u_5\}$ again, before playing on $K_m$, Staller makes another move on $K_m$. Thus a $\gtg^{ssp,4}$-game is played on the cycle. But if Dominator plays on $K_m$ before playing on $\{u_1, u_5\}$ again, at least $\gtg^{sdp (k,\ell)} (C_n)$ or $\gtg^{sdp (k,\ell)} (C_n|^{\ell}\{u_1, u_5\})$ moves are made on the cycle (for some $k \leq l$).
		\end{enumerate}
	\end{enumerate}
	It follows from Lemmas~\ref{lem:cycle15}, \ref{lem:ssp}, \ref{lem:sdp} and~\ref{lem:sdp15} that at least $\frac{2n-1}{3}$ moves are played on $C_n$.	
\end{enumerate}
Hence, $\gtg(G_{n, m}) = \frac{2n-1}{3} + 2$.

Now we describe a strategy for Dominator that shows the inequality $\gtg(G_{n, m}|w_1) \leq \frac{2n-1}{3}$. Dominator starts on $w_1$. After this move the whole $K_m$ (and also vertices $\{u_1, u_5\}$) is totally dominated. Thus all other moves are played on $C_n$, so Dominator follows his optimal strategy on the cycle. Also, Staller plays first on $C_n$. Thus, using total continuation principle, 
$$\gtg(G_{n, m}|w_1) \leq 1 + \gtg'(C_n|\{u_1, u_5\}) \leq 1 + \frac{2n-1}{3} - 1 = \frac{2n-1}{3}\,.$$

As it follows from Lemma~\ref{lem:-2} that $\gtg(G_{n, m}|w_1) \geq \gtg(G_{n,m}) - 2 = \frac{2n-1}{3}$, the proof is complete.

\subsection{Another example}
\label{sec:another}

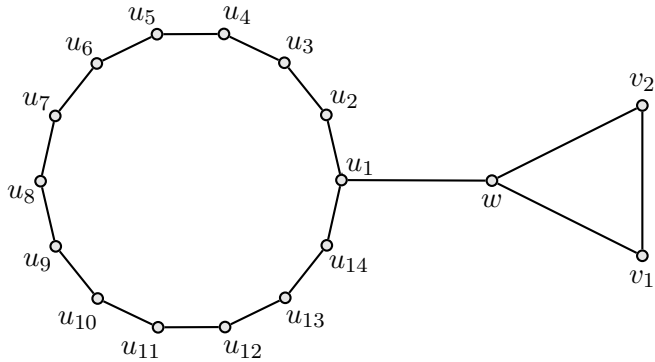
\begin{figure}[!hbt]
 	\begin{center}
 		\begin{tikzpicture}[thick,scale=0.5]
 		
 		\pgfmathtruncatemacro{\N}{14}
 		\pgfmathtruncatemacro{\M}{4}
 		\pgfmathtruncatemacro{\R}{4}
 		\begin{scope}
 		\foreach \x in {1,...,\N}
 		\node[label={-20 + \x*360/\N:$u_{\x}$}] (\x) at (-25.5 + \x*360/\N:\R cm) {};
 		\foreach \x [remember=\x as \lastx (initially 1)] in {1,...,\N,1}
 		\path (\x) edge (\lastx);
 		\end{scope}
 		
 		\node[label=below: {$w$}] (w1) at (8, 0) {};
 		\node[label=below: {$v_1$}] (v1) at (12, -2) {};
 		\node[label=above: {$v_2$}] (v2) at (12, 2) {};
 		
 		\path (1) edge (w1);
 		
 		\path (v1) edge (w1);
 		\path (v2) edge (w1);
 		\path (v2) edge (v1);
 		
 		\end{tikzpicture}
 		\caption{A graph $\widetilde{G}_{14, 3}$.}
 		\label{fig:primer2}
 	\end{center}
 \end{figure}

Another family of graphs with similar properties is the following. Let $\widetilde{G}_{n,m}$, $n, m \geq 3$, be a graph consisting of a cycle $C_n$ with vertices $\U = \{u_1, \ldots, u_n\}$, a complete graph $K_m$ with vertices $w$ and $\V = \{ v_1, \ldots, v_{m-1}\}$ (both with natural edges), and an edge $\{u_1, w\}$ (cf.\ Fig.~\ref{fig:primer2}).

A similar (but simpler) reasoning as above leads to the result $$\gtg(\widetilde{G}_{n,m}) = \frac{2n-1}{3} + 2 \quad \text{ and } \quad \gtg(\widetilde{G}_{n,m}|w) = \frac{2n-1}{3}\,.$$

\section{A generalization}
\label{sec:generalisation}

In Section~\ref{sec:preliminaries} we introduced four new variations of the game total domination number, it would be interesting to find other applications of them. But they are focused to configurations where two vertices, $u_1$ and $u_2$, of a  graph $H$, are connected with two vertices of a complete graph $K_m$. Denote the resulting graph with $G_{H, u_1, u_2, m}$. Following the reasoning of the proof of the Theorem~\ref{thm:minus2}, we can conclude
$$\min_{p,k,l,m}\{\gtg(H|^p\{u_1, u_2\}), \gtg^{ssp, u_1, u_2}(H), \gtg^{sdp(k,l)}(H|^m\{u_1, u_2\})\} + 2 \leq$$
$$\leq \gtg(G_{H, u_1, u_2, m})  \leq \max\{ \gtg(H), \gtg''(H) \} + 2.$$

In the following we consider the case when $H = C_n$. Denote the  graph $G_{H, u_1, u_2, m}$ where $\d(u_1, u_2) = d$ by $G_{n,d,m}$. Notice that $G_{n,m} = G_{n,4,m}$. Using the analogous results as in the preliminaries above, we can state 
\begin{equation}
\label{eqn:meja}
\gtg(C_n|\{u_1, u_2\}) + 2 \leq \gtg(G_{n, \d(u_1, u_2), m}) \leq \gtg(C_n) + 2.
\end{equation}
Therefore we must first determine the value of $\gtg(C_n|\{u_1, u_2\})$ for different choices of vertices $u_1$ and $u_2$. 

\begin{theorem}
	\label{thm:cycles}
	If $n \equiv 2 \mod 6$ and $u_1, u_2 \in V(C_n)$, then
	$$\gtg(C_n|\{u_1, u_2\}) = \begin{cases}
	\gtg(C_n); & \d(u_1, u_2) \mod 6 \equiv 4,\\
	\gtg(C_n) - 1; & \d(u_1, u_2) \mod 6 \in \{1, 3, 5\},\\
	\gtg(C_n) - 2; & \d(u_1, u_2) \mod 6 \in \{0, 2\}.
	\end{cases}$$
\end{theorem}

\proof
We distinguish three cases. Note that $n = 6k + 2$ for some integer $k$.
\begin{enumerate}
	\item $\d(u_1, u_2)$ is odd.\\
		  We first prove the lower bound. Let $A$ denote the set of already totally dominated vertices. By playing as in strategy $S_1$, Staller can always totally dominate only one new vertex. As $u_1, u_2 \in A$ and $\d(u_1, u_2)$ is odd, $(A, A^C)$ can never be a bipartition.
		  
		  Let $m$ denote the number of moves. If $m$ is odd, we have $n \leq 2 + \frac{m-1}{2} + 2 \cdot \frac{m+1}{2}$, thus $m \geq \lceil \frac{2n-5}{3} \rceil = \gtg(C_n) - 1$. If $m$ is even, we get $n \leq 2 + \frac{m}{2} + 2 \cdot \frac{m}{2}$, thus $m \geq \lceil \frac{2n-4}{3} \rceil = \gtg(C_n) - 1$.
		  
		  Consider now the upper bound. Let $v, v' \in V(C_n)$  such that $\d(v, u_1) = \d(v', u_1) = 3$ and $\d(v, u_2) \geq \d(v', u_2)$ (i.e.\ $v$ lies on the longer arc between $u_1$ and $u_2$ and $v'$ on the shorter). Dominator's strategy is to start on $v$ and then reply on the Staller's move such that at least three new vertices become unplayable after these two moves (as in strategy $D_1$). Observe that his reply can be on the same part of the bipartition that Staller played on, unless all vertices in it are unplayable. Notice that two vertices become unplayable after the first Dominator's move and if Staller makes the last move, she also makes two vertices unplayable. 
		  
		  Let $(X, Y)$ be a bipartition of $V(C_n)$ and $u_1 \in X$. Thus $u_2 \in Y$ and $v \in Y$. Notice that $|X| = |Y| = \frac{n}{2} = 3k+1$. Due to the above strategy, Dominator starts playing on $Y$. If Dominator is forced to play first on $X$, he plays the vertex $\widetilde{v} \in V(C_n)$ which is at distance $3$ from $u_2$ and lies on the longer arc between $u_1$ and $u_2$. 
		  
		  Let $l$ denote the number of moves played on $Y$. If $l$ is odd, we have $|X| \geq 2 + 3 \cdot \frac{l-1}{2}$, thus $l \leq 2k-1$. In this case, Staller starts playing on $X$. Let $l'$ denote the number of moves played on $X$. If $l'$ is odd, then $|Y| \geq 3 \cdot \frac{l'-1}{2} + 2$, hence $l' \leq 2k-1$. If $l'$ is even, then $|Y| \geq 3 \cdot \frac{l'}{2}$, so $l' \leq 2k$. In either case, the total number of moves is at most $4k-1 \leq \gtg(C_n) - 1$.
		  
		  If $l$ is odd, we have $|X| \geq 2 + 3 \cdot \frac{l-2}{2} + 2$, so $l \leq 2k$. In this case, Dominator starts on $X$ and repeating the analogous reasoning as for the moves on $Y$, we can conclude that the total number of moves is at most $4k \leq \gtg(C_n) - 1$.
		  
		  So for $\d(u_1, u_2)$ odd, we have $\gtg(C_n|\{u_1, u_2\}) = \gtg(C_n) - 1$.
		  
	\item $\d(u_1, u_2) \equiv 2 \mod 6$.\\
		  Due to strategy $S_1$, we have $\gtg(C_n|\{u_1, u_2\}) \geq \gtg(C_n) - 2$. To show the opposite inequality, we present a suitable strategy for Dominator. Let $(X, Y)$ be the bipartition of $C_n$ and $u_1, u_2 \in X$. Denote $X = \{x_1, \ldots, x_{3k+1}\}, Y = \{y_1, \ldots, y_{3k+1}\}$, where $\d(x_i, x_{i+1}) = \d(y_i, y_{i+1}) = 2$, $x_1 = u_1$ and $y_1, y_2$ are the first two vertices totally dominated on $Y$. Sets $\{x_{3m-1}, x_{3m}, x_{3m+1}\}$ and $\{y_{3m-1}, y_{3m}, y_{3m+1}\}$ are called triplets. Notice that $\{x_1\}$ and $\{y_1\}$ are the only singletons. Let $v \in V(C_n)$ be the vertex at distance $3$ from $u_2$ which lies on the longer arc between $u_1$ and $u_2$. 
		  
		  Dominator's strategy is to start on $v$ and then reply to Staller's moves such that after every two moves, at least one triple becomes totally dominated. Note that this is possible due to the strategy $D_1$. Now consider separately the moves played on $X$ and on $Y$. 
		  
		  After the first move on $Y$, one triple and one singleton is totally dominated. Every next pair of Staller's and Dominator's moves totally dominates at least one new triple. Thus at most $1 + 2 (k-1) = 2k-1$ moves are played on $Y$. 
		  
		  If Staller is the first player to play on $X$, Dominator's first answer should be such that he totally dominates $y_3$ and $y_4$. In this way, the first two moves totally dominate a singleton and a triple. Every next pair of moves totally dominates at least one triple. Hence in this case, the number of moves on $X$ is at most $2 + 2 (k-1) = 2k$. 
		  
		  If Dominator is the first to play on $X$, his move totally dominated only one singleton, and each next pair of moves dominates at least one triple. Hence, the number of moves on $X$ is at most $2k+1$. But this situation only occurs if Dominator is forced to play first on $X$, so there was and even number of moves played on $Y$, thus at most $2k -2$. 
		  
		  In all cases, the total number of moves is at most $4k-1 = \gtg(C_n) - 2$.  So for $\d(u_1, u_2) \equiv 2 \mod 6$, we have $\gtg(C_n|\{u_1, u_2\}) = \gtg(C_n) - 2$.
		  
	\item $\d(u_1, u_2) \equiv 0 \mod 6$.\\
		  Dominator starts on $v$, such that $\d(v, u_2) = 3$ and $v$ lies on the shorter $u_1, u_2$-arc. Now the reasoning is similar as in the previous case.
		  
	\item $\d(u_1, u_2) \equiv 4 \mod 6$.\\
		  It follows from the total continuation principle that $\gtg(C_n|\{u_1, u_2\}) \leq \gtg(C_n)$. Notice that here Dominator cannot choose such first move that exactly one singleton and one triple would be totally dominated after his move. This is probably the reason for a different value of $\gtg(C_n|\{u_1, u_2\})$.
		  
		  Let $(X, Y)$ be the bipartition of $C_n$ (as before) and $u_1, u_2 \in X$. Set $u_1 = x_1$ and $u_2 = x_i$ for some $i$. Staller's strategy is to reply on the same part of the bipartition that Dominator plays. Thus she totally dominates only one new vertex on each move, except if she is the first player to play on $X$. 
		  
		  If Staller is not the first player to play on $X$. \\
		  Thus, Dominator plays first on $X$ and Staller can totally dominate only one new vertex on each of her moves on $X$. Let $l$ be the number of moves Dominator makes on $X$. If he always totally dominated two new vertices, we get $|Y| = 2l + l = 3l$ or $|Y| = 2l + (l-1) = 3l-1$, but neither of those values is congruent to $|Y|$ modulo $3$. Hence, Dominator totally dominates only one new vertex on at least one of his moves. 
		  
		  Let $m$ denote the total number of moves. If $m$ is odd, then $n \leq 2 + \frac{m-1}{2} + 2 \cdot \frac{m+1}{2} - 1$, so $m \geq \frac{2n-1}{3}$. If $m$ is even, then $n \leq 2 + \frac{m}{2} + 2 \cdot \frac{m}{2} - 1$, so $m \geq \frac{2n-1}{3}$.
		  
		  If Staller is the first player to play on $X$.\\
		  This means Dominator started the game on $Y$ and kept playing on $Y$ for as long as possible. Also, he was the last player to play on $Y$. Suppose Dominator totally dominates two new vertices on every move on $Y$. Staller can reply to Dominator's move by totally dominating the remaining vertex in the triple Dominator has just  partially dominated. In the end, exactly one singleton remains undominated on the shorter and one on the longer arc between $u_1$ and $u_2$. But Dominator cannot totally dominate them in one move. Hence, he totally dominates only one new vertex on at least one move. Let $l$ be the number of moves Dominator plays on $Y$. If he totally dominates only one new vertex on just one of his moves, then $|X| - 2 = 2l - 1 + l-1 = 3l-2$, which does not match the size of $X$ modulo $3$. Thus, Dominator dominates just one new vertex on at least two of his moves. 
		  
		  Let $m$ denote the total number of moves. If $m$ is odd, then $n \leq 2 + \frac{m-1}{2} + 1 + 2 \frac{m+1}{2} - 2$, so $m \geq \frac{2n-1}{3}$. If $m$ is even, then $n \leq 2 + \frac{m}{2} + 1 + 2 \frac{m}{2} - 2$, so $m \geq \frac{2n-1}{3}$.
		  
		   So for $\d(u_1, u_2) \equiv 4 \mod 6$, we have $\gtg(C_n|\{u_1, u_2\}) = \gtg(C_n)$. \qed
\end{enumerate}

To complete the study of cycles with two vertices predominated, we also state the following.

\begin{proposition}
	\label{thm:cycles'}
	If $n \equiv 2 \mod 6$ and $u_1, u_2 \in V(C_n)$, then
	$$\gtg'(C_n|\{u_1, u_2\}) = \gtg(C_n) - 1.$$
\end{proposition}

\proof
From the total continuation principle it follows that $\gtg'(C_n|\{u_1, u_2\}) \leq \gtg'(C_n) = \gtg(C_n) - 1$. To prove the lower bound on $\gtg'(C_n|\{u_1, u_2\})$, consider the strategy $S_1$ for Staller. Using similar observations as in the previous proofs, we can see that $\gtg'(C_n|\{u_1, u_2\}) \geq \gtg(C_n) - 1$. 
\qed

Simplifying~\eqref{eqn:meja} yields that for $\d(u_1, u_2) \equiv 4 \mod 6$, it holds $\gtg(G_{n, \d(u_1, u_2), m}) = \gtg(C_n) + 2$. For $\d(u_1, u_2)$ odd, we get $\gtg(C_n) + 1 \leq \gtg(G_{n, \d(u_1, u_2), m}) \leq \gtg(C_n) + 2$, and for $\d(u_1, u_2) \equiv 0 \text{ or } 2 \mod 6$, we have $\gtg(C_n) \leq \gtg(G_{n, \d(u_1, u_2), m}) \leq \gtg(C_n) + 2$. However, computer calculations for $n \in \{8, 14, 20, 26, 32\}$ and $m = 4$ indicate the following.

\begin{problem}
	Is it true that if $n \equiv 2 \mod 6$ and $u_1, u_2 \in V(C_n)$ such that $\d(u_1, u_2) \not \equiv 4 \mod 6$, then 
	$$\gtg(G_{n, \d(u_1, u_2), m}) = \gtg(C_n) + 1 \, ?$$
\end{problem}

The above results give rise to another family of graphs with the property $\gtg(G|v) = \gtg(G) - 2$ for some vertex $v$. Suppose that Dominator's first move on the graph $G_{n, \d(u_1, u_2), m}|w_1$ is on $w_1$, then using Proposition~\ref{thm:cycles'} we get: 
$$\gtg(G_{n, \d(u_1, u_2), m}|w_1) \leq 1 + \gtg'(C_n|\{u_1, u_2\}) = \gtg(C_n).$$ Hence, the family $G_{n, \d(u_1, u_2), m}$ for $\d(u_1, u_2) \equiv 4 \mod 6$, has the desired property (and is in fact a generalization of the family $G_{n,m}$).

\section{Effect of vertex removal on game total domination number}
\label{sec:vertex}

As already mentioned, the effect of vertex removal on the game domination number has been studied in~\cite{EdgeVertexRemoval, book}.  Here we present some analogous results for the game total domination number.

Let $G$ be a graph and $v$ one of its vertices. The game total domination number of the graph  ${G-v}$ cannot be bounded from above by $\gtg(G)$, moreover the difference ${\gtg(G-v) - \gtg(G)}$ can be arbitrarily large. Let $H$ be a graph with $\gtg(H) = k$ and $v \notin V(H)$ a vertex. By connecting $v$ to all other vertices of the graph $H$, we obtain the graph $G$. Clearly, $\gtg(G) = 2$ and $\gtg(G-v) = \gtg(H) = k$. Furthermore, if $H$ is $p$-connected, then $G$ is also $p$-connected. 

But we can bound $\gtg(G-v)$ with $\gtg(G)$ from below. Notice that the bound is weaker compared to the analogous bound for the game domination number. 

\begin{proposition}
	\label{prop:-v}
	If $G$ is a graph and $v \in V(G)$, then 
	$$\gtg(G) \leq \gtg(G-v) + 4.$$
\end{proposition}

\proof
As Dominator can start on $v$, it holds $\gtg(G) \leq 1 + \gtg'(G|N(v))$. Using the imagination strategy we can prove that $\gtg'(G|N(v)) \leq \gtg'((G-v)|N(v)) + 2$. Indeed, it follows from~\cite[Theorem 2.2]{totDom} that $\gtg'(G-v) \leq \gtg(G-v) + 1$. Combining these results and using the total continuation principle yields the desired result.

Hence, it only remains to prove that  $\gtg'(G|N(v)) \leq \gtg'((G-v)|N(v)) + 2$. A total domination game is played on $G|N(v)$, while simultaneously Dominator imagines a Staller-pass game on $(G-v)|N(v)$ (and plays optimally on it). Each move of Staller is copied from the real to the imagined game. This is not possible at most once, exactly when Staller totally dominates only $v$. In this case, Staller passes a move in the imagined game. Dominator replies optimally and copies his move to the real game (this is always legal). Hence, there are at most $\gtg'^{sp}((G-v)|N(v))$ moves made on the imagined game. When it is finished, $v$ might be still undominated in the real game. Thus $\gtg'(G|N(v)) \leq 1 + \gtg'^{sp}((G-v)|N(v))$. 

With the similar reasoning as in~\cite{union} we can derive $\gtg'^{sp}(H) \leq 1 + \gtg'(H)$ for any graph $H$. Thus we have $\gtg'(G|N(v)) \leq 2 + \gtg'((G-v)|N(v))$. 
\qed

A natural question arising from here is whether the bound in Proposition~\ref{prop:-v} is sharp and which differences $\gtg(G) - \gtg(G-v)$ can be realized. We present some partial results on this problem.

\begin{enumerate}
	\item $\gtg(G) - \gtg(G-v) = 0$.
	
	Let $k \in \N$ and $G$ be a graph obtained from $K_{k+2}$, $V(K_{k+2}) = \{ u, v, x_1, \ldots, x_k \}$, by attaching a leaf $y_i$ to $x_i$ for all $i \in [k]$. Notice that both in $G$ and $G-v$, vertices $x_1, \ldots, x_k$ must be played in order to totally dominate all leaves. 
	
	Suppose Dominator starts on $x_1$. If Staller replies on some $x_i$, then the only playable vertices are $\{x_2, \ldots, x_k\} - \{x_i\}$, hence at most $k$ moves are made all together. If Staller replies on $\{ u, v, y_1 \}$, then the only still playable vertices are $\{x_2, \ldots, x_k\}$, thus at most $k+1$ moves are made in total. If Staller replies on some $y_i$, $i \neq 1$, then Dominator replies on $x_i$ and leaves only the vertices $\{x_2, \ldots, x_k\} - \{x_i\}$ playable. Thus again, at most $k+1$ moves are played on the graph.  This strategy for Dominator yields both $\gtg(G) \leq k+1$ and $\gtg(G-v) \leq k+1$. 
	
	Similarly, we observe that $\gtg(G) \geq k+1$ and $\gtg(G-v) \geq k+1$. Hence, $\gtg(G) = \gtg(G-v)$. 
	
	\item $\gtg(G) - \gtg(G-v) = 1$.
	
	It follows from~\cite{totDomCP} that $\gtg(P_n) - \gtg(P_n-v) = 1$ for $n \equiv 0, 1, 2, 4 \mod 6 $ where $v$ is an end-vertex of the path $P_n$.
	
\end{enumerate}

Consider now the Staller-start game. Similarly as in the $D$-game, the value of $\gtg'(G-v)$ cannot be bounded from above by $\gtg'(G)$. But we can determine a lower bound (which is again weaker than for the ordinary domination game).

\begin{proposition}
	\label{prop:-v'}
	Let $G$ be a graph and $v \in V(G)$. Then 
	$$\gtg'(G) \leq \gtg'(G-v) + 4.$$
\end{proposition}

\proof
Using the imagination strategy as in the proof of Proposition~\ref{prop:-v} we can show that $\gtg(G|S) \leq 2 + \gtg((G-v)|S)$. 

If Staller starts on $v$, then we can conclude
$$\gtg'(G) = 1 + \gtg(G|N(v)) \leq 3 + \gtg(G-v | N(v)) \leq 4 + \gtg'(G-v).$$

If Staller starts on a vertex $x$, $x \neq v$, then Dominator's strategy is to reply on $v$. If this move is not legal, we have $N(x) = N(v)$ and without loss of generality, Staller could start on $v$ instead of $x$, thus we have the above situation. But if Dominator can reply on $v$, we have
\begin{eqnarray*}
	\gtg'(G) & = & 1 + \gtg(G | N(x)) \leq 2 + \gtg'(G | N(x) \cup N(v)) \leq \\
	& \leq &  4 + \gtg'(G-v | N(x) \cup N(v)) \leq 4 + \gtg'(G-v),
\end{eqnarray*}

which concludes the proof.
\qed

The question whether the above bound is sharp remains unanswered, but we present some examples of realizable values of $\gtg'(G) - \gtg'(G-v)$. 

\begin{enumerate}
	\item $\gtg'(G) - \gtg'(G-v) = 0$.
	
	Let $k \in \N$ and $G$ be a graph obtained from $K_{k+2}$, $V(K_{k+2}) = \{ u, v, x_1, \ldots, x_k \}$, by attaching a leaf $y_i$ to $x_i$ for all $i \in [k]$. As above, we can prove that $\gtg'(G) = \gtg'(G-v) = k+1$. 
	
	\item $\gtg'(G) - \gtg'(G-v) = 1$.
	
	It follows from~\cite{totDomCP} that $\gtg'(P_n) - \gtg'(P_n-v) = 1$ for $n \equiv 1, 2, 4, 5 \mod 6 $ where $v$ is an end-vertex of the path $P_n$.
	
	\item $\gtg'(G) - \gtg'(G-v) = 2$.
	
	Recall the family of graphs $Z_k$ from~\cite{EdgeVertexRemoval}. Let $Z_0$ be as in Figure~\ref{fig:grafZ}. The graph $Z_k$, $k \geq 1$, is obtained from $Z_0$ by identifying end-vertices of $k$ copies of $P_6$ with the vertex $x$ (cf.~Figure~\ref{fig:grafZk} for $Z_3$). Denote the graph induced by $x$, $u$, $v$, and the two leaves attached to $x$, by $S$. Denote the graph $Z_0 - S$ by $Z$. Observe that $\gamma_t(Z) = \gtg(Z) = \gtg'(Z) = \gtg(Z|z) = \gtg'(Z|z) = 4$ and $\gamma_t(P_6) = 3$.  
	
	\begin{figure}[!hbt]
		\begin{center}
		\begin{minipage}{.5\textwidth}
			\begin{center}
			\begin{tikzpicture}[thick]
			
			\node (1) at (0,1) {};
			\node (2) at (1,1) {};
			\node (3) at (2,1) {};
			\node[label=above: {$z$}] (4) at (3,1) {};
			\node (5) at (0,0) {};
			\node (6) at (1,0) {};
			\node (7) at (2,0) {};
			\node (8) at (3,0) {};
			\node (9) at (2.5,0.5) {};
			\node[label=above: {$x$}] (10) at (4,1) {};
			\node (11) at (3.7,0.5) {};
			\node (12) at (4.3,0.5) {};
			\node[label=-60: {$u$}] (13) at (4,0.5) {};
			\node[label=-60: {$v$}] (14) at (4,0) {};
			
			\draw (1)--(2)--(3)--(4)--(8)--(7)--(6)--(5)--(1);
			\draw (2)--(6); \draw (3)--(7); \draw (4)--(8);
			\draw (9)--(3); \draw (9)--(4)--(10); \draw (9)--(8);
			\draw (11)--(10)--(12); \draw (10)--(13)--(14);
			
			\end{tikzpicture}
			\end{center}
			\caption{A graph $Z_0$.}
			\label{fig:grafZ}
			
		\end{minipage}%
		\begin{minipage}{.5\textwidth}
			\begin{center}
			\begin{tikzpicture}[thick]
			
			\node (1) at (0,1) {};
			\node (2) at (1,1) {};
			\node (3) at (2,1) {};
			\node[label=above: {$z$}] (4) at (3,1) {};
			\node (5) at (0,0) {};
			\node (6) at (1,0) {};
			\node (7) at (2,0) {};
			\node (8) at (3,0) {};
			\node (9) at (2.5,0.5) {};
			\node[label=above: {$x$}] (10) at (4,1) {};
			\node (11) at (3.7,0.5) {};
			\node (12) at (4.3,0.5) {};
			\node[label=-60: {$u$}] (13) at (4,0.5) {};
			\node[label=-60: {$v$}] (14) at (4,0) {};
			\node (15) at (5,1) {};
			\node (16) at (5.5,1) {};
			\node (17) at (6,1) {};
			\node (18) at (6.5,1) {};
			\node (19) at (7,1) {};
			\node (20) at (5,0.7) {};
			\node (21) at (5.5,0.7) {};
			\node (22) at (6,0.7) {};
			\node (23) at (6.5,0.7) {};
			\node (24) at (7,0.7) {};
			\node (25) at (5,0.4) {};
			\node (26) at (5.5,0.4) {};
			\node (27) at (6,0.4) {};
			\node (28) at (6.5,0.4) {};
			\node (29) at (7,0.4) {};
			
			\draw (1)--(2)--(3)--(4)--(8)--(7)--(6)--(5)--(1);
			\draw (2)--(6); \draw (3)--(7); \draw (4)--(8);
			\draw (9)--(3); \draw (9)--(4)--(10); \draw (9)--(8);
			\draw (11)--(10)--(12); \draw (10)--(13)--(14);
			\draw (10)--(15)--(16)--(17)--(18)--(19);
			\draw (10)--(20)--(21)--(22)--(23)--(24);
			\draw (10)--(25)--(26)--(27)--(28)--(29);
			
			\end{tikzpicture}
			\end{center}
			\caption{A graph $Z_3$.}
			\label{fig:grafZk}
		\end{minipage}
		\end{center}
	\end{figure}
	
	We prove that $\gtg'(Z_k) = 3 k + 8$ and $\gtg'(Z_k - v) = 3 k + 6$. Notice that at least four moves are played on $Z$ and at least three moves are played on each path. Consider the following strategy for Staller. She starts the game on $Z_k$ by playing the vertex $v$. If Dominator replies on $u$, then she plays on a neighbor of $x$ on one of the paths. Then she can ensure at least three moves on $S$ and at least four moves on this path. Hence, the total number of moves is at least $3k + 8$. If Dominator replies anywhere else, then she replies optimally on the same subgraph and thus ensures at least four moves on $S$. Hence, $\gtg'(Z_k) \geq 3 k + 8$. 
	
	We now describe a strategy for Dominator to show that $\gtg'(Z_k) \leq 3 k + 8$. If Staller plays on $Z$, then Dominator replies optimally on $Z$. If Staller plays on $S$ or one of the paths, then Dominator replies optimally on the same subgraph. Except immediately after the first move of Staller outside $Z$, when he replies by playing $x$. In this way he ensures that at most $3 k + 8$ moves are played all together.
	
	By applying similar reasoning to the graph $Z_k-v$ we can prove that $\gtg'(Z_k-v) = 3 k + 6$. 
	
\end{enumerate}

An interesting question arising from here is whether there exist graphs $G$ and their vertices $v$ such that $\gtg(G) - \gtg(G-v) \in \{2,3,4\}$ and $\gtg'(G) - \gtg'(G-v) \in \{3,4\}$? And if not, can it be proven in general that for example $\gtg(G) - \gtg(G-v) \leq 2$?

%
%
%
%

\end{document}